\documentclass[12pt]{article}
\usepackage{amsfonts,amssymb,version}

\def\A{{\mathcal A}}
\def\B{{\mathcal B}}

\def\C{{\mathbb C}}
\def\Cc{{\mathcal C}}

\def\E{{\mathcal E}}
\def\EE{{\bf E}}

\def\R{{\mathbb R}}

\def\Z{\mathbb Z}
\def\Zz{{\mathcal Z}}

\def\deg{\mathop{\rm deg}\nolimits}

\def\im{\mathop{\rm im}\nolimits}

\let\dna\downarrow

\let\ov\overline

\let\pa\partial

\makeatletter \def\l@section{\@dottedtocline{1}{0em}{1.2em}} \makeatother

\begin{document}

\centerline{\Large\bf Sign lemma for dimension shifting}

\medskip

\centerline{\bf Nitin Nitsure}


\bigskip

\begin{abstract}
There is a surprising occurrence 
of some minus signs in the isomorphisms produced in the well-known 
technique of dimension shifting in calculating derived functors in 
homological algebra.
We explicitly determine these signs. Getting these signs right 
is important in order to avoid basic contradictions. We illustrate
the lemma by some de Rham cohomology and Chern class considerations
for compact Riemann surfaces. 
\end{abstract}


{\bf \large Statement of the main result}

Let $\A$ be an abelian category with enough 
injectives, and let $F : \A \to \B$ be an additive, left-exact functor
from $\A$ to another abelian category $\B$. For each 
object $M$ of $\A$, we choose an injective resolution
$0\to M \to I^{\bullet}$, and define the value of the 
derived functor $R^iF$ on $M$ to be $H^i(FI^{\bullet})$. 
(To make such a choice for each object, we need to assume some
foundational framework, which is fairly standard so we will omit all
reference to it.) If
$0\to M \to J^{\bullet}$ is any other resolution of $M$, we have
a homomorphism of complexes $f^{\bullet} : J^{\bullet} \to I^{\bullet}$
which is unique up to homotopy, and is identity on 
$M$. If $J^{\bullet}$ is $F$-acyclic, that is, 
if $R^iF$ is zero on $J^k$ for each $i\ge 1$ and $k\ge 0$, 
then for each $n\ge 1$, $f^{\bullet}$ induces an isomorphism
$$c^n = H^n(Ff^{\bullet}) : H^n( FJ^{\bullet}) \to H^n( FI^{\bullet}) 
= R^nFM.$$
We will call $c^n : H^n( FJ^{\bullet}) \to R^nFM$ 
as the {\bf canonical isomorphism} for the 
acyclic resolution $J$ of $M$.

There is another very useful 
isomorphism $d^n : H^n( FJ^{\bullet}) \to R^nFM$ 
for $n\ge 1$, known as the {\bf dimension shifting isomorphism}. 
To define it, 
we begin by breaking-up the resolution $J$ into a sequence 
$\E_1,\ldots,\E_n$ of short exact sequences
\begin{eqnarray*}
\E_1 & = &(0\to Z^0 \to J^0 \to Z^1\to 0), \mbox{ where }Z^0 = M.  \\
\E_2 & = &(0\to Z^1\to J^1 \to Z^2\to 0), \\
&\ldots&\\
\E_n & = &(0\to Z^{n-1}\to J^n \to Z^n\to 0). 
\end{eqnarray*}
As the $J^i$ as $F$-acyclic, 
the corresponding connecting homomorphisms 
are isomorphisms for $p\ge 1$ and $q\ge 1$, which we denote by
$$\pa^p_{\E_q} : R^pFZ^q \to R^{p+1}FZ^{q-1}.$$ 
For $p=0$ and $q\ge 1$, the connecting homomorphism
$\pa^0_{\E_q} : FZ^q \to R^1FZ^{q-1}$ is epic, and induces an isomorphism 
$$\ov{\pa}^0_{\E_q} :  {FZ^q\over \im FJ^{q-1}}  \to 
R^1FZ^{q-1}.$$

For $q=n$, we thus have a sequence of isomorphisms
$$H^n( FJ^{\bullet}) = 
{FZ^n\over \im FJ^n} \stackrel{\ov{\pa}}{\to} R^1F Z^{n-1}
\stackrel{\pa}{\to}  \ldots \stackrel{\pa}{\to} R^nFM$$
The composite of these is an isomorphism
$$d^n : H^n( FJ^{\bullet})\to R^nFM$$
which is by definition the dimension shifting isomorphism for $n\ge 1$.

In this note, we compare the two isomorphisms 
$c^n : H^n( FJ^{\bullet}) \to R^nFM$ and 
$d^n : H^n( FJ^{\bullet})\to R^nFM$ for all $n\ge 1$, 
and we find the following, which is our main result.

\bigskip

{\bf Sign lemma for dimension shifting}
{\it With notation as above, the canonical isomorphism $c^n$ and the 
dimension shifting isomorphism $d^n$ are related by
$$d^n = (-1)^{(n^2+n)/2} \, c^n.$$
}

\bigskip

To prove the above lemma, we need some preliminaries.

\bigskip

\bigskip

{\bf \large Preliminary lemmas} 

Any object $X$ in a category $\Cc$ defines a contravariant functor 
$h_X = Hom_{\Cc}(-,X)$, which we call its {\bf functor of elements}, 
and for any other object $T$, an element $x\in h_X(T)$ will be called a 
{\bf $T$-valued element} of $X$.
When we do not want to explicitly mention $T$, then
such an $x$ will be just called a {\bf valued element} of 
$X$, and by abuse of notation we write it as $x\in X$.
By the Yoneda lemma, any morphism $f: X\to Y$ in $\Cc$ is
determined by its effect on all valued elements of $X$. If $g: X'\to X$
is epic, then (even though not all $T$-valued points 
of $X$ lift to $X'$) the morphism $f: X\to Y$ is
determined by the effect of $g\circ f$ on all valued elements of $X'$.

In what follows, $\A$ will be an abelian category with enough 
injectives, and $F : \A \to \B$ will be an additive, left-exact functor
from $\A$ to another abelian category $\B$.

Let there be given objects $A$ and $B$ in $\A$, together with 
$F$-acyclic resolutions
$0\to A \to J^{\bullet}$ and $0 \to B \to K^{\bullet}$.
Let $c^i_B: H^i(FJ^{\bullet}) \to R^iFA$ and
$c^i_B: H^i(FK^{\bullet}) \to R^iFB$ be the corresponding 
canonical isomorphisms. 
Let there be given a short-exact sequence
$$\E = (0\to A \to C \to B \to 0)$$
and let $\pa_{\E}^i :R^iFB \to R^{i+1}FA$ be the connecting homomorphism. 
Let there be chosen an $F$-acyclic resolution $0\to C\to L^{\bullet}$,
together with a short-exact sequence of complexes
$$\EE = (0 \to J^{\bullet} \to L^{\bullet} \to K^{\bullet} \to 0)$$
such that the following diagram commutes:
$$\begin{array}{ccccc} 
0 \to & A & \to C \to & B & \to 0 \\
      & \dna & \dna& \dna & \\
0 \to & J^{\bullet} & \to L^{\bullet} \to & K^{\bullet} & \to 0
\end{array}$$
Note that by the so called `horse-shoe lemma', such a resolution 
$L^{\bullet}$ always exists. 
As $J^{\bullet}$ is $F$-acyclic, applying $F$ gives a short-exact 
sequence of complexes
$$F\EE = (0 \to FJ^{\bullet} \to FL^{\bullet} \to FK^{\bullet} \to 0)$$
Let $\delta^i_{F\EE} : H^i(FK^{\bullet}) \to H^{i+1}(FJ^{\bullet})$
be the corresponding connecting homomorphism.

\bigskip

{\bf Lemma A.} {\it With notation as above, 
the following diagram commutes.
$$\begin{array}{ccc}
H^i(FK^{\bullet}) & \stackrel{\delta^i_{F\EE}}{\to} & H^{i+1}(FJ^{\bullet})\\
{\scriptstyle c^i}\dna ~~ & & ~~\dna {\scriptstyle c^{i+1}}\\
R^iFA & \stackrel{\pa^i_{\E}}{\to} & R^{i+1}FB
\end{array}$$
}

{\bf Proof.} We leave the proof of Lemma A as an exercise to the reader. 

{\bf Remark.} 
In particular, the homomorphism
$\delta^i_{F\EE} : H^i(FK^{\bullet}) \to H^{i+1}(FJ^{\bullet})$
does not depend on the choice of the resolution
$0\to C\to L^{\bullet}$, or on the short exact sequence
of complexes $\EE$, but only depends
on the given short exact sequence $\E = (0\to A \to C \to B \to 0)$. 
Hence in what follows we will denote it simply by 
$\delta^i_{\E} : H^i(FK^{\bullet})
\to H^{i+1}(FJ^{\bullet})$.

\bigskip

We now return to the situation of our main
result, the Sign Lemma. Recall that we began with an $F$-acyclic 
resolution $0\to M \to J^{\bullet}$ of an object $M$ of $\A$. 
Let $\delta^i_J : J^i\to J^{i+1}$ denote the differentials.
We brake up the resolution into 
short-exact sequences 
$$\E_q = (0 \to Z^{q-1} 
\stackrel{u_{q-1}}{\to} J^{q-1} \stackrel{v_{q-1}} \to Z^q \to 0)$$
for $1\le q\le n$, with $\delta^i = u_{i+1}v_i$. 
For each $i\ge 0$, we get an $F$-acyclic resolution 
$$0 \to Z^i \to K_i^{\bullet}$$
of $Z^i$ defined by
$K_i^p = J^{i+p}$, $\delta^i_K = \delta^{i+p}_J$, and 
$Z^i \to K_i^0$ the `inclusion' homomorphism $u_i: Z^i \to J^i$.

\bigskip

{\bf Lemma B.} {\it For all $n\ge 1$ and $p\ge 1$, 
the following diagrams are commutative.}
$$\begin{array}{ccc}
H^n(FJ^{\bullet}) & \stackrel{-1}{\to} & H^n(FJ^{\bullet}) \\
\|                &                   & \|                \\
{H^0(FK_n^{\bullet})\over \im FJ^{n-1}}
     & \stackrel{\ov{\delta}_{\E_n}}{\to} & H^1(FK_{n-1}^{\bullet})\\
\| & & ~~\dna {\scriptstyle c^1}\\
{FZ^n\over \im FJ^{n-1}}& \stackrel{\ov{\pa}_{\E_n}}{\to} & R^1FZ^{n-1}
\end{array}
\mbox{ and }
\begin{array}{ccc}
H^n(FJ^{\bullet}) & \stackrel{(-1)^{p+1}}{\to} & H^n(FJ^{\bullet}) \\
\|                &                   & \|                \\
H^p(FK_{n-p}^{\bullet})& \stackrel{\delta_{\E_{n-p}}}{\to} 
                            & H^{p+1}(FK_{n-p-1}^{\bullet})\\
{\scriptstyle c^p}\dna ~~ & & ~~\dna {\scriptstyle c^{p+1}}\\
R^pFZ^{n-p}       & \stackrel{\pa_{\E_{n-p}}}{\to} & R^{p+1}FZ^{n-p-1}
\end{array}$$

{\bf Proof} By Lemma A, the homomorphism
$\delta_{\E_q}$ can be computed in terms of any resolution 
$L_i^{\bullet}$ of $J^q$ which fits in a short-exact
sequence $\EE$ of commuting 
resolutions of $0\to Z^q \to J^q \to Z^{q+1}\to 0$, and the
lower squares in the above diagrams commute by Lemma A.
To see that the upper squares commute, 
we construct a particular such resolution 
$$0\to J^i \to L_i^{\bullet}$$
as follows. For any $p$ we put 
$$L_i^p = K_i^p \oplus K_{i+1}^p  = J^{i+p} \oplus J^{i+p+1}.$$
We write valued elements of $L_i^p$ as $(2\times 1)$-column vectors.
With this notation,
the inclusion of $J^i$ into $L_i^0$ is defined in matrix terms by
$$\left(\begin{array}{c}
1_{J^i} \\ \delta^i_J 
\end{array}\right)
: J^i \to J^i \oplus J^{i+1} = L_i^0.$$
The differential 
$$\delta^p_{L_i} : L_i^p \to L_i^{p+1}$$ 
is defined in matrix terms (acting on column vectors) by
$$\delta^p_{L_i} = 
\left(\begin{array}{cc}
\delta^{i+p}_J & (-1)^{p+1} 1_{J^{i+p+1}} \\
0 & \delta^{i+p+1}_J 
\end{array}\right) : J^{i+p} \oplus J^{i+p+1}
\to J^{i+p+1} \oplus J^{i+p+2}.$$
With these definitions, $0\to J^i \to L_i^{\bullet}$ is indeed exact.
Moreover, the following is a commutative 
diagram with exact rows, where the second row is given by
inclusions and projections for the level-wise direct sum 
$L_i^p = K_i^p \oplus K_{i+1}^p  = J^{i+p} \oplus J^{i+p+1}$.
$$\begin{array}{ccccc} 
0 \to & Z^i & \to J^i \to & Z^{i+1} & \to 0 \\
      & \dna & \dna& \dna & \\
0 \to & K_i^{\bullet} & \to L_i^{\bullet} \to & K_{i+1}^{\bullet} & \to 0
\end{array}
$$
Now take $i = n-p -1$ in the above.  
Given any valued element $x\in FZ^n$ which represents
a valued element 
$$\ov{x} \in FZ^n/\im FJ^{n-1} = H^p(FK_{n-p}^{\bullet}),$$ 
the valued element
$$y = 
\left(\begin{array}{c}
0\\ x
\end{array}\right) 
\in  FJ^{n-1}\oplus FJ^n = FK^p_{n-p-1}\oplus FK^p_{n-p} = FL^p_{n-p-1}$$
has the property that under the projection $FL^p_{n-p-1} \to FK_{n-p}^p$,
we have $y\mapsto x$. Now note that as $\delta_{FJ}^nx =0$, we have
$$\left(\begin{array}{cc}
\delta^{n-1}_{FJ} & (-1)^{p+1} 1_{FJ^n} \\
0 & \delta^n_{FJ} 
\end{array}\right)
\left(\begin{array}{c}
0\\ x
\end{array}\right) 
= 
\left(\begin{array}{c}
(-1)^{p+1} x\\  \delta^n_{FJ}x
\end{array}\right)
=
\left(\begin{array}{c}
(-1)^{p+1} x\\ 0
\end{array}\right)
\in FL^{p+1}_{n-p-1}$$
This is the image of $(-1)^{p+1} x \in FK_{n-p-1}^{p+1} = FJ^n$
under the inclusion $FK_{n-p-1}^{p+1} \to FL^{p+1}_{n-p-1}$.
This shows that under the connecting morphism 
$\delta_{F\EE}^p$, the image of $\ov{x}$ is $(-1)^{p+1}\ov{x}$. 
As $FZ^n \to FZ^n/\im FJ^{n-1} = H^p(FK_{n-p}^{\bullet})$ is epic,
this calculation is enough to show that $\delta_{F\EE}^p$
acts as $(-1)^{p+1}$ on all valued elements of $H^p(FK_{n-p}^{\bullet})
= FZ^n/\im FJ^{n-1}$. 
This completes the proof of Lemma B.

\bigskip

\bigskip

{\bf \large Proof of the sign lemma}

By Lemma B, the following squares commute for all $n\ge 1$ and $p\ge 1$.
$$
\begin{array}{ccc}
H^n(FJ^{\bullet}) & \stackrel{-1}{\to} & H^n(FJ^{\bullet}) \\
\| & & ~~\dna {\scriptstyle c^1}\\
{FZ^n\over \im FJ^{n-1}}& \stackrel{\ov{\pa}}{\to} & R^1FZ^{n-1}
\end{array}
\mbox{ and }
\begin{array}{ccc}
H^n(FJ^{\bullet}) & \stackrel{(-1)^{p+1}}{\to} & H^n(FJ^{\bullet}) \\
{\scriptstyle c^p}\dna ~~ & & ~~\dna {\scriptstyle c^{p+1}}\\
R^pFZ^{n-p}  & \stackrel{\pa}{\to} & R^{p+1}FZ^{n-p+1}
\end{array}$$

As $\sum_{p=0}^{n-1} (p+1) = (n^2 +n)/2$, the horizontal composition
of the above diagrams gives a commutative square
$$\begin{array}{ccc}
H^n(FJ^{\bullet}) & \stackrel{(-1)^{(n^2 +n)/2}}{\to} & H^n(FJ^{\bullet}) \\
\|& & ~~\dna {\scriptstyle c^n}\\
{FZ^n\over \im FJ^{n-1}}  & \stackrel{d^n}{\to} & R^nFZ^0
\end{array}$$
where $d^n$ is the composite 
$H^n(FJ^{\bullet}) \stackrel{\ov{\pa}}{\to} R^1F Z^{n-1}
\stackrel{\pa}{\to}  \ldots \stackrel{\pa}{\to} R^nFM$,
which is by definition the dimension-shifting isomorphism.
Thus 
$$d^n = (-1)^{(n^2 +n)/2}\,c^n$$
which completes the proof of the sign lemma.

\bigskip
\bigskip

{\bf \large Illustration: Chern class and de Rham's theorem.}

Let $X$ be a differential manifold (paracompact), and $\A$ the 
category of sheaves of real vector spaces (or complex 
vector spaces) on $X$. Let $\B$
be the category of real (resp. complex) 
vector spaces, and let $F : \A \to \B$
be the global section functor $\Gamma(X,-)$, with derived
functors the sheaf cohomologies 
$H^i(X,-)$. Let $\Cc^i$ be the sheaf of real (resp. complex)
differential $i$-forms, and let $\delta^i: \Cc^i \to \Cc^{i+1}$
be the exterior derivative. As this defines an $F$-acyclic 
resolution $0\to \R_X \to \Cc^{\bullet}$
(resp. $0\to \C_X \to \Cc^{\bullet}$), we get a canonical
isomorphism 
$$c^i: H^i_{dR}(X) = H^i(F\Cc^{\bullet}) \to R^iF(\R_X) = 
H^i(X,\R_X)$$ (a canonical
isomorphism 
$$c^i: H^i_{dR}(X) = H^i(F\Cc^{\bullet}) \to R^iF(\C_X) = 
H^i(X,\C_X)$$
in the complex case), where $H^i_{dR}(X)$ denotes the real (resp. complex)
de Rham cohomology of $X$.

Some authors (for example [G-H]) prove  
de Rham's theorem by identifying de Rham cohomology and
sheaf cohomology by the dimension shifting
isomorphism 
$d^i : H^i_{dR}(X) \to H^i(X,\C_X)$.  
By the sign lemma, this is $(-1)^{(i^2+i)/2}$-times the canonical 
isomorphism. Omission of this sign can lead to
sign mistakes and confusion later on. 

One such confusion, which we
describe next to end this note, occurs in the basic calculation
of Chern classes of line bundles on Riemann surfaces.
Let $\Cc^*$ be the 
multiplicative sheaf non-vanishing complex valued smooth functions
on $X$.
Let $\exp: \Cc\to \Cc^*$ the map defined at the level
of local sections by $f\mapsto e^{2\pi i f}$. This map is 
surjective at the level of germs, and defines a short exact sequence
of sheaves
$$0\to \Z_X \to \Cc \stackrel{\exp}{\to} \Cc^* \to 0$$
called as exponential sequence. Any complex 
line bundle on $X$ defines an element $(L) \in H^1(X, \Cc^*)$,
whose image 
$$c_1(L) = \pa(L) \in  H^2(X,\Z_X)$$
under the connecting homomorphism 
$\pa : H^1(X, \Cc^*)\to H^2(X,\Z_X)$
is the first Chern class of $L$.

For a compact Riemann surface $X$,
let $\eta_X \in H^2(X,\Z_X)$ denote the positive generator. 
For a complex line bundle $L$ on a 
compact Riemann surface it can be directly calculated that 
$$\pa(L) = \deg(L)\,\eta_X \in  H^2(X,\Z_X)$$
where $\deg(L)$ is the degree of $L$, which is positive for ample
line bundles. Hence we get the relation
$c_1(L) = \deg(L)\,\eta_X$ relating first Chern class and degree.

We have a commutative diagram with
exact rows
$$\begin{array}{ccccc}
0 \to & \Z_X  & \to \Cc \to  & \Cc^*  & \to 0 \\
      & \dna  &      \|        & \dna     &        \\
0 \to & \C_X & \to \Cc \to  &   \Zz^1  & \to 0
\end{array}$$
where $\Zz^1$ is the sheaf of closed $1$-forms, and 
$\Cc^*_X \to \Zz^1$ is defined by
$f\mapsto df/2\pi i f$.
This gives a commutative diagram
$$\begin{array}{ccc}
H^1(X,\Cc^*) & \stackrel{\pa}{\to} & H^2(X,\Z_X) \\
\dna         && \dna \\      
H^1(X, \Zz^1)& \stackrel{\pa}{\to} &H^2(X,\C_X)
\end{array}$$
in which the bottom row is an isomorphism.
Thus, any complex line bundle $L$ defined by transition functions
$(g_{a,b}) \in H^1(X,\Cc^*)$
defines a class $c(L) = (dg_{a,b}/ 2\pi ig_{a,b})\in  H^1(X, \Zz^1)$
whose image $\pa(c(L)) \in H^2(X,\C_X)$ is $c_1(L)$.
We have connecting isomorphisms
$$H^2_{dR}(X)  = {H^0(X,\Zz^2)\over \im H^0(X, \Cc^1)}
\stackrel{\pa}{\to} H^1(X,\Zz^1) \stackrel{\pa}{\to}
H^2(X,\C_X)$$
whose composite is the dimension shifting isomorphism
$d^2 : H^2_{dR}(X) \to H^2(X,\C_X)$.

A simple calculation (see for example [Na]) shows that 
if $X$ is a compact Riemann surface 
and if $\alpha_X \in H^2_{dR}(X)$ is the positive integral generator
(means $\int_X \alpha_X =1$), then 
$$c(L) = -\deg(L)\, \pa(\alpha_X).$$
This is consistent with the sign lemma, by which 
$d^2 : H^2_{dR}(X) \to H^2(X,\C_X)$ is $(-1)$-times
the canonical isomorphism, so that $\pa\circ \pa (\alpha_X) = -\eta_X$,
and so $\pa(c(L)) = \deg(L)\, \eta_X$.

Not taking the sign $(-1)$ into account will lead to a paradox at this point. 
Any attempt to resolve it by trying
to define the first Chern class as $-\pa(L) \in H^2(X,\Z_X)$ will 
in turn be contradicted by a direct calculation of
$\pa : H^1(X, \Cc^*)\to H^2(X,\Z_X)$. 
The author of this note actually got into this contradiction,
which led to this work which in particular resolves it.

\bigskip
\bigskip

{\bf References}

[G-H] Griffiths, P. and Harris, J. : {\it Algebraic Geometry} 
Wiley-Interscience Pub., 1978.

[Na] Narasimhan, M.S. `Vector bundles on compact Riemann surfaces',
in {\it Complex analysis and its applications} Vol. 3, 
IAEA Vienna, 1976.

\bigskip

\bigskip

{\footnotesize 

\hfill School of Mathematics

\hfill Tata Institute of Fundamental Research

\hfill Homi Bhabha Road

\hfill Mumbai 400 005, India

\hfill e-mail: nitsure@math.tifr.res.in

\centerline{15 June 2007}

}

\end{document}